\documentclass{siamart220329}

\usepackage{arxiv}

\usepackage[utf8]{inputenc} 
\usepackage[T1]{fontenc}    
\usepackage{url}            
\usepackage{booktabs}       
\usepackage{amsfonts}       
\usepackage{nicefrac}       
\usepackage{microtype}      
\usepackage{lipsum}
\usepackage{graphicx}
\graphicspath{ {./images/} }

\usepackage{rotating} 

\usepackage{ragged2e}

\usepackage{scalerel}
\usepackage{tikz}
\usetikzlibrary{svg.path}

\definecolor{orcidlogocol}{HTML}{A6CE39}
\tikzset{
  orcidlogo/.pic={
    \fill[orcidlogocol] svg{M256,128c0,70.7-57.3,128-128,128C57.3,256,0,198.7,0,128C0,57.3,57.3,0,128,0C198.7,0,256,57.3,256,128z};
    \fill[white] svg{M86.3,186.2H70.9V79.1h15.4v48.4V186.2z}
                 svg{M108.9,79.1h41.6c39.6,0,57,28.3,57,53.6c0,27.5-21.5,53.6-56.8,53.6h-41.8V79.1z M124.3,172.4h24.5c34.9,0,42.9-26.5,42.9-39.7c0-21.5-13.7-39.7-43.7-39.7h-23.7V172.4z}
                 svg{M88.7,56.8c0,5.5-4.5,10.1-10.1,10.1c-5.6,0-10.1-4.6-10.1-10.1c0-5.6,4.5-10.1,10.1-10.1C84.2,46.7,88.7,51.3,88.7,56.8z};
  }
}

\newcommand\orcidicon[1]{\href{https://orcid.org/#1}{\mbox{\scalerel*{
\begin{tikzpicture}[yscale=-1,transform shape]
\pic{orcidlogo};
\end{tikzpicture}
}{|}}}}

\usepackage{bm} 
\newcommand{\bx}{{\bf x}}

\newcommand{\bX}{{\bf X}}
\newcommand{\bY}{{\bf Y}}

\newcommand{\bA}{{\bf A}}

\newcommand{\bB}{{\bf B}}
\newcommand{\bC}{{\bf C}}
\newcommand{\bE}{{\bf E}}

\newcommand{\bI}{{\bf I}}
\newcommand{\bK}{{\bf K}}

\newcommand{\bM}{{\bf M}}
\newcommand{\bMetax}{{\bM_{\bmeta \bx}}}

\newcommand{\bS}{{\bf S}}

\newcommand{\bV}{{\bf V}}

\newcommand{\bQ}{{\bf Q}}

\newcommand{\bmxi}{{\bm{\xi}}}

\newcommand{\bmgam}{{\bm{\gamma}}}

\newcommand{\gammaSTc}{{\gamma^\text{C-LIM}}}

\newcommand{\gammanLIM}{{\gamma^\text{C-nLIM}}}
\newcommand{\bmeta}{{\bm{\eta}}}

\newcommand{\Alim}{{A^{\text{LIM}}}}
\newcommand{\AWnLIM}{{A^{\text{W-nLIM}}}}
\newcommand{\ACnLIM}{{A^{\text{C-nLIM}}}}

\newcommand{\ASTw}{{A^{\text{W-LIM}}}}
\newcommand{\ASTc}{{A^{\text{C-LIM}}}}

\newcommand{\BWnLIM}{{B^{\text{W-nLIM}}}}
\newcommand{\BCnLIM}{{B^{\text{C-nLIM}}}}
\newcommand{\Clim}{{C^{\text{LIM}}}}
\newcommand{\CWnLIM}{{C^{\text{W-nLIM}}}}
\newcommand{\CCnLIM}{{C^{\text{C-nLIM}}}}
\newcommand{\Qlim}{{Q^{\text{LIM}}}}
\newcommand{\QWnLIM}{{Q^{\text{W-nLIM}}}}
\newcommand{\QCnLIM}{{Q^{\text{C-nLIM}}}}

\newcommand{\QSTw}{{Q^{\text{W-LIM}}}}
\newcommand{\QSTc}{{Q^{\text{C-LIM}}}}

\newcommand{\LFP}{\textbf{L}_{\textit{FP}}}
\newcommand{\pst}{P_{\textit{st}}}

\newcommand{\Eobs}{{E^{\text{obs}}}}

\newcommand{\Kobs}{{K^{\text{obs}}}}

\newcommand{\Mobs}{{M^{\text{obs}}}}
\newcommand{\Sobs}{{S^{\text{obs}}}}

\usepackage{algorithm,algpseudocode}
\algblock{Input}{EndInput}
\algnotext{EndInput}
\algblock{Output}{EndOutput}
\algnotext{EndOutput}

\usepackage{lipsum}
\usepackage{amsfonts}
\usepackage{graphicx}
\usepackage{epstopdf}
\ifpdf
  \DeclareGraphicsExtensions{.eps,.pdf,.png,.jpg}
\else
  \DeclareGraphicsExtensions{.eps}
\fi

\newcommand{\R}{\mathbb{R}}

\newsiamremark{remark}{Remark}
\newsiamremark{hypothesis}{Hypothesis}
\crefname{hypothesis}{Hypothesis}{Hypotheses}
\newsiamthm{claim}{Claim}

\usepackage{amsmath} 
\usepackage{mathtools} 
\usepackage{booktabs,graphicx,makecell,siunitx,eqparbox,threeparttable}
\usepackage{subcaption}
\usepackage{adjustbox}
\usepackage{multirow}
\usepackage{lscape} 
\usepackage{afterpage}
\usepackage{graphicx} 
\usepackage{array}

\usepackage{tikz}
\usetikzlibrary{arrows.meta}
\tikzset{mynode/.style={draw, very thick, circle, minimum size=0.8cm}, myarrow/.style={very thick}}

\usepackage{graphicx}

\DeclareMathOperator*{\argmin}{arg\,min} 

\usepackage{verbatim}

\usepackage{hyperref} 

\title{ Beyond Gaussian Assumptions: A Nonlinear Generalization of Linear Inverse Modeling }

\author{
    Justin Lien \\
    Mathematical Institute \\
    Tohoku University\\
    Sendai, Japan \\
    \texttt{lien.justin.t8@dc.tohoku.ac.jp} \\
    \And
    Hiroyasu Ando \\
    Advanced Institute for Materials Research\\
    Tohoku University\\
    Sendai, Japan \\
}

\begin{document}
\maketitle
\begin{abstract}


The Linear Inverse Model (LIM) is a class of data-driven methods that construct approximate linear stochastic models to represent complex observational data. 
The stochastic forcing can be modeled using either Gaussian white noise or Ornstein-Uhlenbeck colored noise; the corresponding models are called White-LIM and Colored-LIM, respectively. 
Although LIMs are widely applied in climate sciences, they inherently approximate observed distributions as Gaussian, limiting their ability to capture asymmetries.

In this study, we extend LIMs to incorporate nonlinear dynamics, introducing White-nLIM and Colored-nLIM which allow for a more flexible and accurate representation of complex dynamics from observations. 
The proposed methods not only account for the nonlinear nature of the underlying system but also effectively capture the skewness of the observed distribution.
Moreover, we apply these methods to a lower-dimensional representation of ENSO and demonstrate that both White-nLIM and Colored-nLIM successfully capture its nonlinear characteristic.

\end{abstract}







\section{ Introduction }

The most fundamental linear inverse model (LIM), referred to as White-LIM in this article, is widely used to study local dynamical behavior near an equilibrium point \cite{Penland1989,Penland1996,Penland1993,Penland1994}. 
More specifically, it extracts the local dynamics from observational data by constructing an approximate linear stochastic differential equation (SDE) driven by Gaussian white noise. 
LIM has been extensively applied to study the El Nino-Southern Oscillation (ENSO), effectively capturing key characteristics of ocean-atmosphere coupled variability \cite{Aiken2013,Alexander2008,Lou2020}. 
However, fundamental limitations remain, particularly regarding the predictability and the asymmetry of ENSO events \cite{DelSole1996,DelSole2000,Martinez2017,Shin2021}.

Recently, stochastic forcing has been extended to include noise persistence using Ornstein-Uhlenbeck (OU) colored noise, leading to Colored-LIM \cite{Lien2024}. 
By accounting for noise memory, this approach offers a more realistic representation of an unknown complex system with improved predictability.
However, as a Gaussian model, it still fails to capture ENSO asymmetry.
Meanwhile, a regression-based method, XRO, has been proposed to model ENSO by incorporating second-order dynamics and noise persistence \cite{Zhao2024}. 
XRO effectively captures ENSO asymmetry and achieves forecast skill comparable to state-of-the-art machine learning techniques, highlighting the importance of higher-order dynamics and persistent stochastic forcing in ENSO modeling \cite{Ham2019}.

In this study, we extend both White-LIM and Colored-LIM to a non-linear framework by incorporating quadratic dynamics, and call them White-nLIM and Colored-nLIM. 
We show that regardless of the noise type, the higher-order local dynamics of the quadratic stochastic dynamical system is characterized by the local behavior of its correlation functions up to 4th order. 
When applied to observational data, this framework enables a more accurate and realistic representation of local dynamics near an equilibrium point or a strange attractor. 
In addition, it captures the skewness of the observed probability distribution, addressing key limitations of traditional LIM approaches.

Before proceeding, we clarify the notation and setup used in this study. Boldface denotes mathematical objects in the \text{forward} formulation, while regular font represents estimated quantities from the \text{inverse} method.
Vectors are assumed to be column vectors.
An $a$-tensor $\bX$ is an element in $\R^{\overbrace{ n\times\cdots \times n}^{a \text{ times}}}$.
For an $a$-tensor $\bX$ and a $b$-tensor $\bY$, their tensor multiplication produces an $(a+b-2m)$-tensor $\bX \times_{m} \bY$, by summing over the last $m$ indices of $\bX$ and the first two indices of $\bY$.
For example, given a $3$-tensor $\bX$ and a $4$-tensor $\bY$, their product $\bX \times_{2} \bY$ is a $3$-tensor given by
\[ [\bX \times_{2} \bY]_{ijk} =\sum_{p,\,q} \bX_{i p q} \bY_{p q j k}. \]
The symmetrization operator $\text{Sym}$ for an $a$-tensor $\bX$ is defined by 
\[ [\text{Sym}(\bX)]_{i_1 i_2 i_3\dots i_a} = \bX_{i_1 i_2 i_3\dots i_a} + \bX_{i_2 i_1 i_3\dots i_a}. \]
Throughout this article, stationarity is assumed implicitly, and the SDE can be interpreted in either Ito or Stratonovich sense due to additive noise.

\section{ Review of Linear Inverse Modelings } 
\label{Chap:LIM}
Linear inverse modelling describes the $n$-dimensional observational dataset $\{ x(t): t = 0,\dots,N\Delta t \} \subset \R^n$ that satisfies zero-mean condition (i.e., $\sum x(t) = 0$) by producing an approximate linear stochastic system 
\begin{align} \label{Eq:Approx-form}
    \frac{d}{dt} x = \Alim x + \sqrt{2\Qlim} \cdot \text{noise}.
\end{align}
In this formulation, $\Alim \in \R^{n\times n}$ is a dissipative matrix (i.e., all eigenvalues have negative real parts) that describes the first-order deterministic dynamics of the underlying system, and $\Qlim \in \R^{n\times n}$ is a positive definite matrix that represents the spatial-temporal structure of the stochastic forcing modeled by either white or colored noise. 

\subsection{ White-LIM }

The White-LIM considers an $n$-dimensional Markov system of the form
\begin{align} \label{Eq:ST-White-LIM}
    \frac{d}{dt} \bx = \bA\bx + \sqrt{2\bQ} \bmxi.
\end{align}
The stochastic forcing $\bmxi$ is the normalized Gaussian white noise with zeros mean and unit variance, i.e., $\langle \bmxi(\cdot+\tau)\bmxi(\cdot)^T \rangle = \delta(\tau)\;\bI$,
where the bracket, $\delta$, and $\tau$ denote the expectation, Dirac delta function, and the lag variable, respectively \cite{Penland1989,Penland1995}.
Eq. (\ref{Eq:ST-White-LIM}) is an OU process whose correlation function $\bK$ given by $\bK(\tau) \coloneqq \langle \bx(\cdot+\tau) \bx(\cdot)^T \rangle$ is an exponential function of the form
\begin{align} \label{Eq:ST-K}
    \bK(\tau) = e^{\tau \bA} \bK(0),
\end{align} 
for $\tau \ge 0$, and hence,
\begin{align} \label{Eq:ST-White-Dyn}
    \bK'(0) = \bA\bK(0).
\end{align}
The balanced equation describing the dynamical equilibrium of the deterministic damping and stochastic excitation is given by the fluctuation-dissipation relation (FDR):
\begin{align} \label{Eq:ST-FDR}
    0 = \text{Sym}\left( \, \bA\bK(0) + \bQ \, \right).
\end{align}
Therefore, given that $\bK$ is known, the dynamical and stochastic matrices can be reconstructed using Eqs. (\ref{Eq:ST-White-Dyn}) and (\ref{Eq:ST-FDR}).

In practical applications, given an observational dataset $\{x(t)\}$, the White-LIM makes the Markov assumption, describing the data by Eq. (\ref{Eq:ST-White-LIM}) and estimating the first-order dynamics $\ASTw$ by using Eq. (\ref{Eq:ST-White-Dyn}) and stochastic matrix $\QSTw$ by using Eq. (\ref{Eq:ST-FDR}) with the observed correlation function $\Kobs$ computed by 
\begin{equation} \label{Eq:Kobs}
    \Kobs(\tau) = \frac{\sum_{t=0}^{(N-k)\Delta t} x(t+\tau)x(t)^T}{N-k+1},
\end{equation}
where $\tau = k\Delta t$.

\subsection{ Colored-LIM }

In the Colored-LIM framework, the governing equation is given by 
\begin{align} \label{Eq:ST-Colored-LIM}
    \frac{d}{dt} \bx = \bA\bx + \sqrt{2\bQ} \bmeta.
\end{align} 
The stochastic forcing is modeled by the OU colored noise $\bmeta$, realized by the state steady of 
\begin{equation} \label{Eq:Colored-noise}
    \frac{d}{dt}\bmeta = -\frac{1}{\bmgam} \bmeta + \frac{1}{\bmgam} \bmxi,
\end{equation} 
where $\bmgam$ is the noise correlation time.
The temporal structure of $\bmeta$ is given by 
\begin{equation} 
    \langle \bmeta(\cdot+\tau) \bmeta(\cdot)^T \rangle = \frac{1}{2\bmgam} e^{ \frac{-\left\vert \tau \right\vert}{\bmgam} }\bI,
\end{equation} 
and hence, $\bmgam$ can be physically interpreted as the memory of noise, the persistence of noise, or the characteristic time for the past dependence of noise.

A straightforward computation shows that the covariance $\bK_{\bmeta\bx}$ between the colored noise and state variable satisfies
\begin{equation} 
    \bK_{\bmeta\bx} \coloneqq \langle \bmeta(\cdot)\bx(\cdot)^T \rangle = \sqrt{\frac{\bQ}{2}}\bB^T,
\end{equation}
where $\bB=(\bI-\bmgam\bA)^{-1}$, and the correlation function $\bK = \langle \bx(\cdot+\tau) \bx(\cdot)^T \rangle$ of the state variable is given by
\begin{equation} \label{Eq:ST-Colored-K}
    \bK(\tau) 
    = e^{\tau\bA}\bC + e^{\tau \bA} \int_0^\tau e^{-s(\bA+\frac{1}{\bmgam})} ds \, \bQ\bB^T.
\end{equation}

It has been shown that the second- and third-order derivatives of $\bK$ can be used to reconstruct the linear dynamics \cite{Lien2024}.
However, the second-order derivative alone is sufficient. 
More preciously, $\bK''$ satisfies
\begin{equation} \label{Eq:Colored-Dynamics}
    \bK''(0) = \bA\bK'(0) - \frac{1}{\bmgam}\big(\bK'(0) - \bA\bK(0)\big),
\end{equation}
which is independent of $\bQ$. 
The balanced equation is then described by
\begin{equation} \label{Eq:Colored-FDR}
    0 = \text{Sym} \left( \bA \bK(0) + \bQ\bB^T \right),
\end{equation}
referred to as the Colored-FDR.
Therefore, $\bK$ and $\bmgam$ characterize $\bA$ and $\bQ$ via Eqs. (\ref{Eq:Colored-Dynamics}) and (\ref{Eq:Colored-FDR}).

Given an observation dataset, the Colored-LIM assumes stationarity and models the environmental noise with an OU colored noise.
For a prescribed $\gammaSTc$, it estimates the dynamics $\ASTc$ and stochastic forcing $\QSTc$ by solving Eqs (\ref{Eq:Colored-Dynamics}) and (\ref{Eq:Colored-FDR}) with the observed correlation $\Kobs$ computed by Eq. (\ref{Eq:Kobs}).
When not prescribed, $\gammaSTc$ is determined by Colored-LIM via the $\gamma$-selection algorithm
\begin{align} \label{Eq:tau-selection}
    \gammaSTc = \argmin_{\gamma > 0} \left\Vert K\left(A(\gamma),Q(\gamma),\gamma\right) - \Kobs \right\Vert,
\end{align}
where $A(\tau)$ and $Q(\tau)$ are the solutions of Eqs. (\ref{Eq:Colored-Dynamics}) and (\ref{Eq:Colored-FDR}) with the running varaible $\gamma$, and $\left\Vert \cdot \right\Vert$ is the function norm over the interval $[0,l]$ with $l$ being a parameter controlling the minimization window.
Then we have $\ASTc = A(\gammaSTc)$ and $\QSTc = Q(\gammaSTc)$.

\section{Mathematical Foundations of Non-Linear Inverse Modelings}

The non-linear inverse models extend LIMs to a non-linear framework by incorporating quadratic and constant terms, describing the evolution of the observations $\{x(t)\}$ that may not satisfy the zero-mean condition by
\begin{align} \label{Eq:nLIM-Approx-form}
    \frac{d}{dt} x_i = \sum_{j,\,k} B^\text{nLIM}_{ijk} x_j x_k+\sum_j A^\text{nLIM}_{ij} x_j +C^\text{nLIM}_i + \sum_j\left( \sqrt{2 Q^\text{nLIM}} \right)_{ij} \cdot \text{noise}_j,
\end{align}
where $B^\text{nLIM}$ is a $3$-tensor with $B^\text{nLIM}_{ijk} = 0$ if $j>k$.
The random forcing can be either modeled by Gaussian white noise or OU colored noise. 
Again, we first focus on the forward formulation before addressing the inverse problem.

\subsection{White-nLIM}

In the Wthie-nLIM framework, the governing equation is given by
\begin{align} \label{Eq:White-nLIM}
    \frac{d}{dt} \bx_i = \sum_{j,\,k}\bB_{ijk}\bx_j \bx_k+\sum_j\bA_{ij}\bx_j +\bC_i + \sum_j\left( \sqrt{2\bQ} \right)_{ij} \bmxi_j.
\end{align}
where $\bB_{ijk} = 0$ if $j>k$, $\bQ$ is positive definite, and $\bmxi$ is the normalized white noise.

Due to the quadratic term, the dynamics $\bB$, $\bA$, and $\bC$ are not determined solely by the correlation function $\bK_{ij}(\tau) = \langle \bx_i(\cdot + \tau) \bx_j(\cdot) \rangle$.
Instead, the expectation $\bE_i = \langle \bx_i(\cdot) \rangle$, and the higher-order correlation functions $\bM_{ijk}(\tau) = \langle\bx_i(\cdot + \tau) \bx_j(\cdot)\bx_k(\cdot)\rangle$ and $\bS_{ijkw}(\tau) = \langle\bx_i(\cdot + \tau) \bx_j(\cdot)\bx_k(\cdot)\bx_w(\cdot)\rangle$ come into play.
Despite the lack of explicit forms, these correlation functions satisfy
\begin{align}
    0 &= \bB \times_2 \bK(0) + \bA\times_1\bE  + \bC \label{Eq:White-nLIM-E} \\
    \bK'(0) &= \bB \times_{2}\bM(0) + \bA\times_1\bK(0) +\bC\times_0\bE \label{Eq:White-nLIM-K'} \\
    \bM'(0) &= \bB\times_{2}\bS(0) + \bA\times_1\bM(0) + \bC\times_0\bK(0). \label{Eq:White-nLIM-M'}
\end{align}
Moreover, the balanced equation is given by
\begin{align} \label{Eq:FDR-White-nLIM}
    0 = \text{Sym}\left(\, \bB\times_2\bM(0) + \bA\times_1\bK(0) + \bC\times_0\bE^T+ \bQ \, \right).
\end{align}
Therefore, by viewing the expectation as a degenerate first-order correlation function, we may say that the dynamics are characterized by the correlation functions up to 4th order, and so is the stochastic matrix.

Given observations, the White-nLIM computes the expectation $\Eobs$ and the observed correlation functions $\Kobs$, $\Mobs$, and $\Sobs$ in the same manner as in Eq. (\ref{Eq:Kobs}) as well as their first-order derivatives, and uses the linear equations determined by Eqs. (\ref{Eq:White-nLIM-E}) to (\ref{Eq:White-nLIM-M'}) to estimate the dynamics $\BWnLIM$, $\AWnLIM$, and $\CWnLIM$. 
Then the stochastic matrix $\QWnLIM$ is simply solved by Eq. (\ref{Eq:FDR-White-nLIM}). 

\subsection{Colored-nLIM}

The Colored-nLIM models the stochastic forcing with OU colored noise $\bmeta$, and the governing equation reads
\begin{align} \label{Eq:Colored-nLIM}
    \frac{d}{dt} \bx_i = \sum_{j,\,k}\bB_{ijk}\bx_j \bx_k+\sum_j\bA_{ij}\bx_j + \bC_i + \sum_j\left( \sqrt{2\bQ} \right)_{ij} \bmeta_j
\end{align}
where $\bmeta$ is the colored noise with noise correlation time $\bmgam$.

As in Colored-LIM, the covariance between colored noise and observable is non-trivial.
Let $\bK_{\bmeta\bx} \in \R^{n\times n}$ and $\bM_{\bmeta\bx}\in \R^{n\times n \times n}$ denote the noise-observable correlations, given by
$[\bK_{\bmeta\bx}]_{ij} = \langle\bmeta_i(\cdot) \bx_j(\cdot)\rangle$, and
$[\bM_{\bmeta\bx}]_{ijk} = \langle\bmeta_i(\cdot) \bx_j(\cdot)\bx_k(\cdot)\rangle$, respectively.
Similar to White-nLIM, the correlation functions do not admit explicit forms, but satisfy
\begin{align}
    0 &= \bB \times_2 \bK(0) + \bA\times_1\bE  + \bC \label{Eq:Colored-nLIM-E} \\
    \bK'(0) &= \bB\times_2\bM(0) + \bA\times_1\bK(0) + \sqrt{2\bQ}\times_1\bK_{\bmeta \bx} +\bC\times_0\bE \label{Eq:Colored-nLIM-K'} \\
    \bK''(0) &= \bB\times_2\text{Sym}(\bM'(0)) + \bA\times_1\bK'(0) - \frac{1}{\bmgam} \sqrt{2\bQ}\times_1\bK_{\bmeta \bx} \\
    \bM'(0) &= \bB\times_2\bS(0) +  \bA\times_1\bM(0) + \sqrt{2\bQ} \times_1 \bMetax + \bC \times_0 \bK(0)  \label{Eq:Colored-nLIM-M'} \\
    \bM''(0) &= \bB\times_2 \text{Sym}(\bS'(0)) + \bA \times_1 \bM'(0) - \frac{1}{\bmgam} \sqrt{2\bQ}\times_1\bMetax \label{Eq:Colored-nLIM-M''}
\end{align}
Moreover, we have
\begin{align} \label{Eq:Colored-nLIM-Q}
    0=\bM_{\bmeta\bx} \times_2\bB + \bK_{\bmeta\bx} \times_1\bA^T - \frac{1}{\bmgam}\bK_{\bmeta\bx} + \frac{1}{2\bmgam}\sqrt{2\bQ},
\end{align}
and the balance equation reads
\begin{align} \label{Eq:Colored-nLIM-FDR}
    0 = \text{Sym}\left(\, \bB\times_2\bM(0) + \bA \times_1\bK(0) + \bC\times_0\bE  + \sqrt{2\bQ} \times_1 \bK_{\bmeta\bx} \, \right).
\end{align}
Consequently, the dynamics $\bB$, $\bA$, and $\bC$ as well as the stochastic matrix $\bQ$ are functions of $\bE$, $\bK$, $\bM$, $\bS$, and $\bmgam$.
    
Given observations $\{x(t)\}$ and a prescribed noise correlation time $\gammanLIM$, the estimated parameters $\BCnLIM$, $\ACnLIM$, and $\CCnLIM$ are obtained by solving the linear system determined by Eqs. (\ref{Eq:Colored-nLIM-E}) to (\ref{Eq:Colored-nLIM-M''}). 
The estimated stochastic matrix $Q_0$ then is solved by Eqs. (\ref{Eq:Colored-nLIM-K'}), (\ref{Eq:Colored-nLIM-M'}), and (\ref{Eq:Colored-nLIM-Q}). 
If $\{x(t)\}$ is indeed sampled from Eq. (\ref{Eq:Colored-nLIM}), then $Q_0$ indeed provides a fair estimate.

However, in practice, several factors can lead to discrepancies.
These include limited sampling, deviations from the ideal quadratic system, or numerical errors in computing derivatives.
As a result, the reconstructed covariance of the approximate system determined by $\BCnLIM$, $\ACnLIM$, $\CCnLIM$, $Q_0$, and $\gammanLIM$, may not fully align with the observed covariance $\Kobs$ and $\Mobs$.
Therefore, we adopt the minimization scheme as follows:
\begin{align} \label{Eq:Colored-nLIM-minQ}
    \QCnLIM = \argmin_{Q} \left( \left\Vert R\right\Vert_F + \left\Vert K(0)-\Kobs(0) \right\Vert_F + \left\Vert M(0)-\Mobs(0) \right\Vert_F \right)
\end{align}
with the initial guess $Q_0$, where $\left\Vert\cdot\right\Vert_F$ denotes the Frobenius norm and the residual matrix $R$ is defined as
\begin{align}
    R = \text{Sym} \left( \BCnLIM\times_2 \Mobs(0) + \ACnLIM\times_1\Kobs(0) + \CCnLIM \times_0 \Eobs + \sqrt{2\QCnLIM} \times_1 K_{\eta x}  \right).
\end{align}
The matrices $K$, $M$, and $K_{\eta x}$ are the correlation functions and covariance between the state variable and colored noise computed by a sufficiently long realization of the system determined by $\BCnLIM$, $\ACnLIM$, $\CCnLIM$, $\gammanLIM$, and the running variable $Q$.

\subsection{Stationarity}

In practice, the quadratic system (\ref{Eq:nLIM-Approx-form}) determined by White-nLIM or Colored-nLIM may not admit a stationary distribution. 
Nonetheless, one of the primary goals of nLIM is to analyze the local behavior of the underlying system without requiring global stability of the approximate system.

To maintain theoretical rigor and justify such an approximate system, we introduce an artificial boundary to Eq. (\ref{Eq:nLIM-Approx-form}) to bound and reflect the escaping trajectories whenever necessary.
Specifically, this is achieved by adding a dissipative force $\bV$ defined by $\bV(\Vert x\Vert) = -w(\Vert x - x_0\Vert-r)\cdot (x-x_0)$ where $r$ and $x_0$ specify the domain of interest and the function $w$ is given by
\begin{align}
    w(z) = 
    \begin{cases}
        0, &  z < 0\\
        e^{mz}\cdot e^{-z^{-1}}, & z \ge 0
    \end{cases}
    \, ,
\end{align}
for some $m \gg 1$.
The presence of $\bV$ ensures the existence of a stationary distribution of the modified system, and so long as $\langle w(x) \rangle \ll 1$, the influence of the artificial wall is minimal and negligible in practice.
As a result, all key equations in White-nLIM and Colored-nLIM are still effective in constructing the approximate system in practical applications.

\subsection{Relationships among models}

From the perspective of forward formulation, even if the system parameters $\bA$ and $\bQ$ remain unchanged in the LIM framework, the correlation function of Eq. (\ref{Eq:ST-Colored-LIM}) does not converge to that of Eq. (\ref{Eq:ST-White-LIM}) in the white-noise limit $\bmgam \to 0$ \cite{Lien2024}.
However, when applied to the same observational data, the Colored-LIM results do converge to the White-LIM results. 
This can be seen by multiplying $\bmgam$ on both sides of Eq. (\ref{Eq:Colored-Dynamics}) and then taking the white-noise limit, which reduces Eq. (\ref{Eq:Colored-Dynamics}) to Eq. (\ref{Eq:ST-White-Dyn}) and the Colored-FDR converges to FDR.
This suggests that in the inverse problem, the parameter $\gamma$ is a measure of how dominantly the concavity of the observed correlation functions $\Kobs$ and $\Mobs$ influences the approximate system. 
The same reasoning applies to nLIMs. 

The LIM framework can also be viewed as a special case of the nLIM framework. 
In the forward formulation, if the quadratic and constant terms vanish, the key equations of nLIMs reduce to those of LIMs.
Furthermore, for white-noise models, when applied to the same observational data, the stochastic matrix depends only on the first-order derivative of $\Kobs$, making it independent of the presence of nonlinear and constant terms.
In contrast, for colored-noise models, the stochastic matrices $\QSTc$ and $\QCnLIM$ explicitly depend on higher-order observed correlation functions and therefore differ.
Finally, when data is sampled sufficiently close to a stable equilibrium point of the underlying system, LIMs are known to approximate the Jacobian matrix \cite{Lien2024}. 
Under this condition, the nLIM framework yields an estimated linear term that closely aligns with the result of the linear models.

\section{ Validation } 
In this section, we assess the effectiveness of White-nLIM and Colored-nLIM.
All numerical derivatives in this study are computed using finite forward difference schemes.
To efficiently represent the $3$-tensor $\bB\in\R^{n\times n\times n}$ of quadratic dynamics, we identify it with an ordinary $n$-by-$n(n+1)/2$ matrix whose $i$th-row is the transpose of the vectorization of the upper triangular part of $[\bB_i]_{jk} = \bB_{ijk}$. That is,
\begin{align}
    \bB = 
    \begin{bmatrix}
        \bB_{111} & \bB_{112} & \cdots & \bB_{11n} & \bB_{122}  & \cdots & \bB_{1nn} \\
        \vdots & & & \ddots& & & \vdots \\
        \bB_{n11} & \bB_{n12} & & \cdots& & & \bB_{nnn}
    \end{bmatrix}
    .
\end{align}

It is well-known that a minimum dimension for a differentiable autonomous system to exhibit chaotic behavior is $3$ \cite{Sandri1996}. 
Therefore, we start with a $2$-dimensional example and then the chaotic Lorenz 63 system to validate the effectiveness of White-nLIM and Colored-nLIM to non-chaotic and chaotic systems.

\subsection{ A 2-dimensional example }

We consider the following non-linear systems:
\begin{align} \label{Eq:Example}
    \bB = 
    \begin{bmatrix}
        0 & -1 & 0 \\
        1 & 0 & 0
    \end{bmatrix}
    \text{, }
    \bA = 
    \begin{bmatrix}
        -1  & 2 \\
        -1 & -2
    \end{bmatrix}
    ,
    \bC = 
    \begin{bmatrix}
        0.5 \\
        0
    \end{bmatrix}
    \text{, and }
    \bQ = 
    \begin{bmatrix}
        1  & 0.5 \\
        0.5 & 1
    \end{bmatrix}
    ,
\end{align}
subject to white-noise and colored-noise stochastic forcing.
In the case of colored noise, we set $\bmgam = 0.5$ and assume it is a known parameter.

Following \cite{Lien2024}, we generate 100 realizations of the system determined by Eq. (\ref{Eq:Example}) using the Euler-Maruyama method with an integration time step $dt = 0.001$ and a total time span of $T_f = 1000$ and $5000$.
The data is then subsampled so that the observations $\{x(t)\}$ have equal sampling interval $\Delta t = 10\,dt$. 
Finally, we apply nLIMs to estimate the dynamics and stochastic matrix.
To evaluate the performance, we compute the relative error $e$ and absolute error $E$ using the Frobenius norm for a given tensor $X_1$ and its reference $X_0$:
\[ e_{(X_0,X_1)} = \frac{\left\Vert X_0 - X_1 \right\Vert_F}{\left\Vert X_0 \right\Vert_F}, \]
and 
\[ E_{(X_0,X_1)} = \left\Vert X_0 - X_1 \right\Vert_F. \]

Both models effectively estimate the system parameters in this setup, with 
relative errors summarized in Table \ref{Table:Abs_Error}. 
As discussed in \cite{Lien2024}, the performance of LIMs depends on (a) whether the observed correlation functions, especially $\Kobs$ and $\Mobs$, accurately reflect the true correlation functions and (b) whether the numerical derivatives are accurately computed. 
To investigate the source of relative errors, we numerically compute the true correlation functions $\bK$ and $\bM$ as well as the correlation functions $K^\text{nLIM}$ and $M^\text{nLIM}$ determined by nLIM, and compare them with the observed ones using absolute errors.

The absolute errors summarized in Table \ref{Table:Abs_Error} suggest that nLIMs can accurately fit the $\Kobs$ and $\Mobs$, without referring to $\bK$ and $\bM$.
This supports the idea that the primary source of relative errors is limited sampling. 
If an infinite number of observations were available, then the observed correlation functions approach the true ones, and the numerical errors would be solely due to the finite difference scheme used in computing numerical derivatives. 
This result is consistent with \cite{Lien2024}.

\begin{table}[htbp] 

\sisetup{
  table-space-text-post=\%,    
  table-align-text-post=false, 
}
\begin{centering}
\footnotesize
\caption{The median of the numerical errors for the system determined by (\ref{Eq:Example}). The true correlation functions $\bK$ and $\bM$ as well as those determined by nLIMs $K^\text{nLIM}$ and $M^\text{nLIM}$ are computed by synthetic data of length $10\, T_f$.}
\label{Table:Abs_Error}
\begin{threeparttable}

\begin{tabular}{ l S[table-format=-2.5] S[table-format=-2.5] S[table-format=-2.5] S[table-format=-2.5] }
    \toprule
    \multicolumn{1}{c}{} & \multicolumn{2}{c}{White-nLIM} & \multicolumn{2}{c}{Colored-nLIM}\\
    \multicolumn{1}{c}{} & \multicolumn{1}{c}{$T_f = 1000$} & \multicolumn{1}{c}{$T_f = 5000$} & \multicolumn{1}{c}{$T_f = 1000$} & \multicolumn{1}{c}{$T_f = 5000$} \\
    \midrule
    \makecell{ $e_{(\bB,\,B^\text{nLIM})}$ } & 9.90\si{\percent} & 5.07\% & 17.77\% & 7.65\% \\ 
    \makecell{ $e_{(\bA,\,A^\text{nLIM})}$ } & 3.82\% & 2.14\% & 5.35\% & 2.54\% \\ 
    \makecell{ $e_{(\bC,\,C^\text{nLIM})}$ } & 14.63\% & 7.53\% & 15.42\% & 6.61\% \\ 
    \makecell{ $e_{(\bQ,\,Q^\text{nLIM})}$ } & 0.53\% & 0.36\% & 3.98\% & 2.65\% \\ 
    \midrule
    \makecell{ $E_{(\Kobs,K^\text{nLIM})}$ } & 0.3734 & 0.1827 & 0.4415 & 0.1747  \\ 
    \makecell{ $E_{(\Kobs,\,\bK)}$ } & 1.1328 & 0.4487 & 1.0450 & 0.4201 \\ 
    \midrule
    \makecell{ $E_{(\Mobs,M^\text{nLIM})}$ } & 0.9178 & 0.5054 & 1.0645 & 0.3825 \\ 
    \makecell{ $E_{(\Mobs,\,\bM)}$ } & 3.2730 & 1.2425 & 2.5797 & 1.0502 \\ 
    \bottomrule
\end{tabular} 
\begin{tablenotes}
    \item 
\end{tablenotes}
\end{threeparttable}

\end{centering}
\end{table}

\subsection{ Stochastic Lorenz 63 system }

The Lorenz 63 system consists of three coupled deterministic ordinary differential equations that exhibit chaotic behavior \cite{Lorenz63}.
We explore its stochastic counterpart, introducing stochastic forcing with colored noise defined by $\bmgam = 0.5$ and $\bQ = \bI$. 
More precisely, the system is governed by:
\begin{align}
    \frac{d}{dt}\bx = F(\bx) + \sqrt{2\bQ} \,\bmeta,
    \text{ where } 
    \begin{cases}
        F_1(\bx) &= \bm{\sigma} (\bx_2-\bx_1) \\
        F_2(\bx) &= \bx_1(\bm{\rho}-\bx_3)-\bx_2 \\
        F_3(\bx) &= \bx_1\bx_2 - \bm{\beta}\bx_3 \\
    \end{cases}
\end{align}
with parameters $\bm{\sigma} = 10$, $\bm{\rho} = 28$, and $\bm{\beta} = 8/3$.
The goal here is to estimate the original deterministic parameters using observations from the stochastic system.

Given its chaotic nature, we use a 4th-order Runge-Kutta method to integrate the deterministic part and the Euler-Maruyama method to update the colored noise \cite{LeVeque2007,sarkka_solin_2019}.
The system is evolved using a time step $dt = 0.001$ over a total duration $T_f = 10000$, with ($\Delta t = 10\,dt$) and without ($\Delta t = dt$) subsampling to generate the synthetic data $\{x(t)\}$.

Table \ref{Table:Lorenz63} summarizes the estimated parameters by an application of Colored-nLIM with and without the restriction
\begin{align} \label{Eq:Restriction}
    \begin{cases}
        B^\text{LIM}_{213},B^\text{LIM}_{312} \ne 0 \\
        B^\text{LIM}_{ijk}= 0 \text{, otherwise}
    \end{cases} 
    \,,\,
        \begin{cases}
        A^\text{LIM}_{11} = -A^\text{LIM}_{12} \\
        A^\text{LIM}_{21}, A^\text{LIM}_{22}, A^\text{LIM}_{33} \ne 0 \\
        A^\text{LIM}_{ij} = 0 \text{, otherwise}
    \end{cases}
    , 
    \text{ and } \Clim = 0.
\end{align}
The estimated parameters become inaccurate if the synthetic data is generated with subsampling and without restriction (\ref{Eq:Restriction}).
To investigate the source of this error, we construct a sufficiently long high-resolution auxiliary dataset to numerically approximate the true correlation functions $\bK$ and $\bM$.
Our numerical computation shows that coarse subsampling causes up to 15\% relative errors in the numerical derivatives as they are sensitive to the data resolution probably due to the chaotic nature of the Lorentz 63 system.
Moreover, the linear system that determines the parameters in the Colored-nLIM framework may not be perfectly well-conditioned. 
To mitigate the numerical errors, we may increase the sampling frequency or reduce the number of estimated parameters, as shown in Table \ref{Table:Lorenz63}.

\begin{table}[htbp] 
\footnotesize
\begin{centering}
\footnotesize
\caption{The relative errors of the Lorenz 63 system for a single realization of $T_f = 10000$.}
\label{Table:Lorenz63}
\begin{threeparttable}
\footnotesize
\begin{tabular}{ l l }
    \toprule
    \multicolumn{1}{c}{} & \multicolumn{1}{c}{Estimated parameters} \\
    \midrule
    \makecell{ No restriction } & \makecell{ $\BCnLIM = $ \( \begin{bmatrix*}[r]
    -0.0001 & 0.0001 & -0.0043 & -0.0001 & -0.0007 & -0.0000 \\
     0.0007 & -0.0146 & -0.9631 & 0.0120 & -0.0279 & 0.0015 \\
    -0.0059 & 0.9993 & -0.0003 & 0.0000 & -0.0003 & -0.0027
    \end{bmatrix*} \), \\
    $\ACnLIM = $ 
    \( \begin{bmatrix*}[r]
    -9.8338 &  9.9668 &  0.0014 \\
    26.6609 & -0.1024 & -0.0055 \\
     0.0169 & -0.0002 & -2.5756 \\
    \end{bmatrix*} \), and 
    $\CCnLIM = $
    \( \begin{bmatrix*}[r]
    -0.0144 \\
     0.0166 \\
    -0.2893 \\
    \end{bmatrix*} \) }
    \\ 
    \midrule
    \makecell{ With restriction } & \makecell{ $\sigma^\text{C-nLIM} = 10.0191$, $\rho^\text{C-nLIM} = 29.2046$, $\beta^\text{C-nLIM} = -2.6808$, \\ $B^\text{C-nLIM}_{213} = -1.1735$, $B^\text{C-nLIM}_{312} = 0.9907$, $A^\text{C-nLIM}_{22} = -1.0251$ } \\ 
    \midrule
    \makecell{ No restriction, subsampling } & \makecell{ $\BCnLIM = $ \( \begin{bmatrix*}[r]
    -0.0000 & -0.0003 & -0.0340 & 0.0002 & -0.0105 &  0.0000 \\
    -0.0011 &  0.0007 & -0.6751 & 0.0000 & -0.2344 &  0.0003 \\
    -0.0003 &  0.7268 & -0.0019 & 0.1871 &  0.0011 & -0.0093
    \end{bmatrix*} \), \\
    $\ACnLIM = $ 
    \( \begin{bmatrix*}[r]
    -8.6831 &  9.8083 &  0.0002 \\
    16.4602 &  6.2191 & -0.0155 \\
     0.0721 & -0.0451 & -1.8962
    \end{bmatrix*} \), and 
    $\CCnLIM = $
    \( \begin{bmatrix*}[r]
    -0.0067 \\
    -0.0024 \\
    -0.2475
    \end{bmatrix*} \) }
    \\ 
    \midrule
    \makecell{ With restriction, subsampling } & \makecell{ $\sigma^\text{C-nLIM} = 10.1048$, $\rho^\text{C-nLIM} = 26.7796$, $\beta^\text{C-nLIM} = -2.8460$, \\ $B^\text{C-nLIM}_{213} = -1.1864$, $B^\text{C-nLIM}_{312} = 0.9013$, $A^\text{C-nLIM}_{22} = -0.9742$ }  \\ 
    \bottomrule
\end{tabular} 
\begin{tablenotes}
    \item 
\end{tablenotes}
\end{threeparttable}

\end{centering}
\end{table} 

Overall, our analysis demonstrates that Colored-nLIM can effectively deal with systems with chaotic dynamics, provided that the sampling resolution is sufficiently high and numerical derivatives are accurately computed. 
While Colored-nLIM is theoretically well-suited for chaotic systems, its practical reliability depends on the quality of observations, which may not always be robust.
Finally, a similar analysis conducted for White-nLIM, in which colored noise is replaced by white noise, leads to the same conclusions, so we do not pursue further discussion on White-nLIM.

\section{ Modeling El Nino 3.4 index using LIMs }

ENSO is a dominant mode of climate variability in the tropical Pacific, characterized by warm and cold phases called El Nino and La Nina, respectively, influencing global weather patterns \cite{Lorenzo2015,Wang2017}. 
Sea surface temperature (SST) is a primary indicator of ENSO, driven by oceanic and atmospheric processes, while the thermocline represents the heat exchange between surface and deeper waters \cite{Wen1999}. 
Previous studies have suggested that the interaction between SST and thermocline forms the foundation of ENSO evolution, making them essential variables in simple ENSO models that capture key feedback mechanisms governing ENSO dynamics \cite{Burgers2005,Jin2007}

In this study, we revisit the simple ENSO model consisting of the area-averaged SST and thermocline, defined as the 20$^\circ$C isotherm
depth (D20). 
The data is collected from the ORAS5 reanalysis from 1979 to 2021.
SST is averaged over 170W-120W and 5S-5N, while D20 is averaged over 120E-80W and 5S-5N, following \cite{Lien2024}.
The linear trends and climatological means of SST and D20 are subtracted to remove long-term trends. 
Since LIMs and nLIMs do not account for seasonality, we take z-scores for each calendar month to eliminate seasonal dependence \cite{Kido2023,Lien2024,Shin2021}. 
Finally, a 3-month sliding average is applied to smooth the data \cite{Newman2011}.
The preprocessed dataset is denoted by $\{x(t)\} \subset \R^2$, where SST and D20 are indexed by $1$ and $2$, respectively.

\begin{table}[htbp] 
\begin{centering}
\footnotesize
\caption{Model results}
\label{Table:Model-Results}
\begin{threeparttable}

\begin{tabular}{ l l l l l l}
    \toprule
    \multicolumn{1}{c}{Model} & \multicolumn{1}{c}{$B$} & \multicolumn{1}{c}{$A$} & \multicolumn{1}{c}{$C$} & \multicolumn{1}{c}{$Q$} & \multicolumn{1}{c}{$\gamma$} \\
    \midrule
    \makecell{ White-LIM } & \multicolumn{1}{c}{-} & 
    \( \begin{bmatrix*}[r]
         -0.7770 & 1.6444 \\
         -1.6018 & 0.0392
    \end{bmatrix*} \) & \multicolumn{1}{c}{-} & 
    \( \begin{bmatrix*}[r]
        0.3660 & 0.0709 \\
        0.0709 & 0.3612
    \end{bmatrix*} \) & \multicolumn{1}{c}{-}
    \\ 
    \midrule
    \makecell{ Colored-LIM } & \multicolumn{1}{c}{-} & 
    \( \begin{bmatrix*}[r]
        -1.1646 & 1.6850 \\
        -1.6091 & -0.2849
    \end{bmatrix*} \) & \multicolumn{1}{c}{-} & 
    \( \begin{bmatrix*}[r]
        0.8047 & 0.1560 \\
        0.1560 & 0.7321 
    \end{bmatrix*} \)  
    & \multicolumn{1}{c}{0.10}
    \\
    \midrule
    \makecell{ White-nLIM } & 
    \( \begin{bmatrix*}[r]
        0.0973 & 0.5784 & 0.0488 \\
        -0.3181 & 0.4038 & -0.1346
    \end{bmatrix*} \) & 
    \( \begin{bmatrix*}[r]
        -0.7700 & 1.7691 \\
        -1.3908 & -0.1407
    \end{bmatrix*} \) &
    \( \begin{bmatrix*}[r]
        -0.2906 \\
        0.3518
    \end{bmatrix*} \) & 
    \( \begin{bmatrix*}[r]
        0.3660 & 0.0709 \\
        0.0709 & 0.3612
    \end{bmatrix*} \)  &
    \multicolumn{1}{c}{-}
    \\ 
    \midrule
    \makecell{ Colored-nLIM } & 
     \( \begin{bmatrix*}[r]
        0.1184 & 0.7847 & 0.0235 \\       
        -0.4439 & 0.5046 & -0.3066
    \end{bmatrix*} \) & 
    \( \begin{bmatrix*}[r]
        -1.1025 & 1.8082 \\
        -1.3182 & -0.6284
    \end{bmatrix*} \) &
    \( \begin{bmatrix*}[r]
        -0.3380 \\
         0.6243
    \end{bmatrix*} \) & 
    \( \begin{bmatrix*}[r]
        0.7074 & 0.1975 \\
        0.1975 & 0.6705
    \end{bmatrix*} \)  & \multicolumn{1}{c}{0.08} \\
    \midrule
    \bottomrule
\end{tabular} 
\begin{tablenotes}
    \item The unit of time is in months, and the unit of state variable is dimensionless.
\end{tablenotes}
\end{threeparttable}

\end{centering}
\end{table} 

Table \ref{Table:Model-Results} summarizes the model outputs from LIMs and nLIMs, where the noise correlation time $\gammaSTc$ in Colored-LIM is determined by Eq. (\ref{Eq:tau-selection}) with minimization window $l = 12$ months, and $\gammanLIM$ is prescribed for Colored-nLIM.
Regardless of the noise persistence, the results that $B^\text{LIM}_{111}$ and $B^\text{LIM}_{112} > B^\text{LIM}_{122}$ are consistent with the XRO studies, where the quadratic terms are linked to deterministic nonlinear ocean advection and atmospheric nonlinearity implicitly through the nonlinear SST-wind stress feedback \cite{Zhao2024}.
On the other hand, the other quadratic terms $B^\text{LIM}_{2jk}$, assumed to be zero in the XRO studies, are non-zero and comparable to $B^\text{LIM}_{112}$, suggesting the presence of other non-linear processes contributing to subsurface dynamics.

\begin{figure}[htbp]
\centering 
\includegraphics{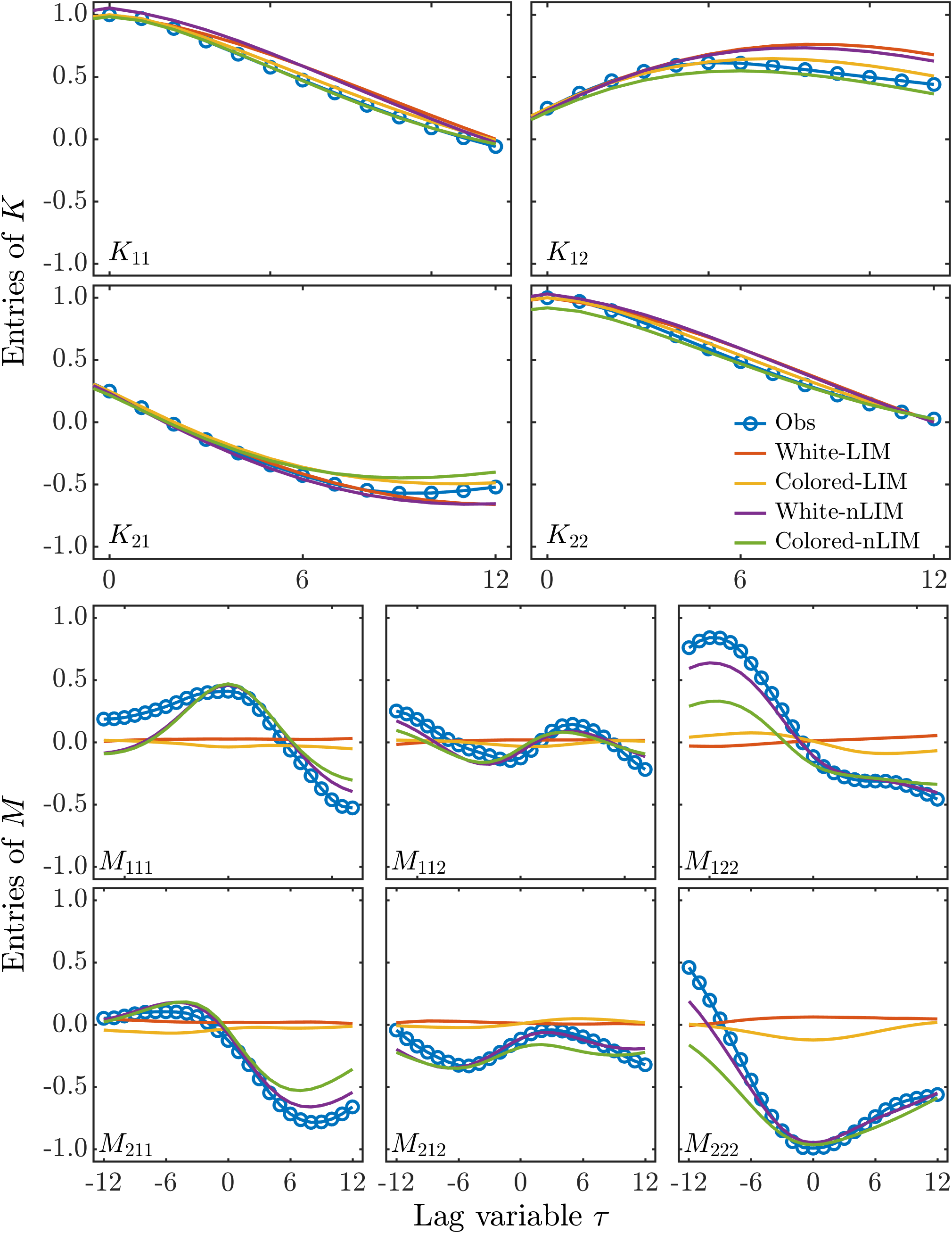}
    \caption{The entries of the observed correlation function $\Kobs$ and $\Mobs$ of the preprocessed SST and D20 anomalies symbolized as 1, and 2, respectively, along with the correlation functions determined by White-LIM, Colored-LIM, White-nLIM, and Colored-nLIM. For LIMs, the correlation functions are obtained from Eqs. (\ref{Eq:ST-K}) and (\ref{Eq:ST-Colored-K}), while for nLIMs, they are determined by a 1000-year realization of the approximate systems.}
    \label{Fig:Corr}
\end{figure}

Figure \ref{Fig:Corr} shows the observed correlation functions $\Kobs$ and $\Mobs$ along with those determined by the approximate systems. 
Since all models learn from $\Kobs$ by design, they inherently respect the observed covariance. 
However, nLIMs go further by incorporating information from $\Mobs$, allowing them to capture its variation in the lag variables.
In contrast, the linear models do not capture the higher-order observed correlation functions.

Compared to the white-noise models, the colored-noise models provide a better fit to $\Kobs$, especially near the origin.
This improvement arises because the correlation functions of the colored-noise models are local parabolas in the diagonal entries, as noted in \cite{DelSole1996,Lien2024}, as well as indicated by Eqs. (\ref{Eq:Colored-nLIM-K'}), and (\ref{Eq:Colored-nLIM-FDR}), due to the non-trivial $\gamma$. 
On the other hand, white-noise models produce correlation functions with cusps, a feature that is less pronounced in $\Kobs$. 

\begin{figure}[htbp]
\centering 
\includegraphics{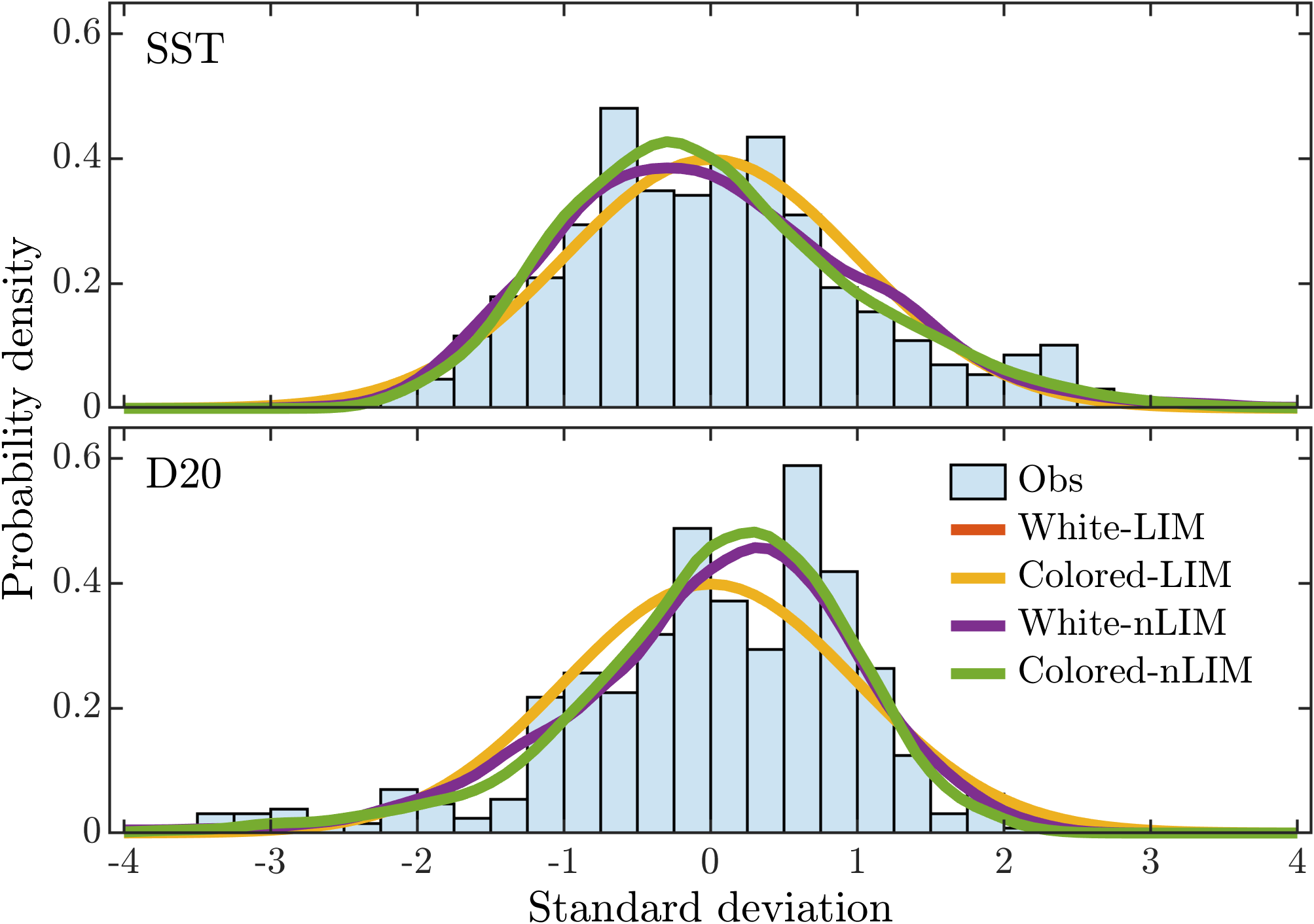}
    \caption{The observed distributions of preprocessed SST and D20, along with the probability distributions of the approximate systems determined by White-LIM, Colored-LIM, White-nLIM, and Colored-nLIM. White-LIM and Colored-LIM produce Gaussian distributions whose covariances are the diagonal entries of $\Kobs(0)$ and they are indistinguishable. For White-nLIM and Colored-nLIM, the distributions are determined by a 1000-year realization of the approximate systems. }
    \label{Fig:Dist}
\end{figure}

Figure \ref{Fig:Dist} shows the observed distribution and the probability distributions of the approximate systems determined by LIMs and nLIMs, and Table \ref{Table:Ws-distance} shows their Wasserstein $p$-distance.
The results indicate that the skewed nLIM distributions align more closely with the observed distributions, in contrast to the symmetric Gaussian distributions produced by linear models. 
Given this skewed characteristic, nLIMs offer the potential to better represent and facilitate the study of asymmetric El Nino and La Nina events.

\begin{table}[htbp] 

\sisetup{detect-weight=true,detect-inline-weight=math}
\begin{centering}
\footnotesize
\caption{Wasserstein $p$-distances between the observed distribution and the probability distribution determined by LIMs and nLIMs.}
\label{Table:Ws-distance}
\begin{threeparttable}

\begin{tabular}{ l S[table-format=2.5] S[table-format=2.5] S[table-format=2.5] S[table-format=2.5] }
    \toprule
    \multicolumn{1}{c}{} & \multicolumn{2}{c}{$p = 1$} & \multicolumn{2}{c}{$p = 2$} \\
    \multicolumn{1}{c}{Model} & \multicolumn{1}{c}{SST} & \multicolumn{1}{c}{D20} & \multicolumn{1}{c}{SST} & \multicolumn{1}{c}{D20} \\
    \midrule
    \makecell{ White-LIM } & 0.0949 & 0.1826 & 0.1354 & 0.2397 \\
    \midrule
    \makecell{ Colored-LIM } & 0.0949 & 0.1826 & 0.1354 & 0.2397 \\
    \midrule
    \makecell{ White-nLIM } & 0.0713 & 0.0753 & 0.1187 & \bfseries 0.1380 \\
    \midrule
    \makecell{ Colored-nLIM } & \bfseries 0.0555 & \bfseries 0.0670 & \bfseries 0.0924 & 0.1488 \\
    \midrule
    \bottomrule
\end{tabular} 
\begin{tablenotes}
    \item The Wasserstein $p$-distance is unitless.
\end{tablenotes}
\end{threeparttable}

\end{centering}
\end{table}

\section{Concluding Remarks}

In this study, we have extended linear inverse modeling to a nonlinear framework by incorporating quadratic nonlinear dynamics and constant terms, as well as modeling random forcing as either memoryless Gaussian white noise or persistent OU colored noise, referred to as White-nLIM and Colored-nLIM, respectively.
The numerical experiments confirm that despite being derivative-based algorithms, nLIMs are effective when the expectation, observed correlation functions, and their derivatives accurately capture the characteristics of the underlying system.

In principle, the nLIM framework is applicable to chaotic systems, but the reliability depends on the quality of observations. 
We have demonstrated that the estimated parameters can be incorrect if numerical derivatives are not sufficiently accurate or if the parameter estimation process is not well-conditioned.
This suggests that imposing constraints on the estimation may be necessary to improve numerical stability. 

Applied to real-world time-series data for ENSO, both White-nLIM and Colored-nLIM faithfully capture the observed correlations and the skewed distributions of SST and D20.
This is a notable improvement over conventional linear models, which assume symmetric Gaussian distributions and therefore cannot reproduce this skewness. 

Although the stationary XRO shares conceptual similarities to Colored-nLIM, the two approaches differ in key aspects. 
XRO utilizes regression techniques to produce the approximate system \cite{Zhao2024}.
As the colored noise is correlated with the state variables, standard regression methods may lead to biased estimates, requiring more advanced regression approaches.
In contrast, Colored-nLIM is primarily based on the local behaviors of correlation functions, so a robust algorithm may be needed to compute numerical derivatives. 

Furthermore, a critical distinction is that XRO does not include a constant term in the deterministic dynamics, whereas Colored-nLIM does. 
This distinction is particularly important in ENSO studies, where observational data typically represent anomalies of climate variables, meaning their mean states are inherently constrained to be zero.
To ensure consistency with the zero-mean condition, Colored-nLIM (and White-nLIM) explicitly incorporates the constant term in the governing equation. 

In the Colored-nLIM framework, the noise correlation time $\gammanLIM$ is treated as a hyperparameter. 
The primary challenge of a $\gamma$-selection algorithm in the same manner as in Eq. (\ref{Eq:tau-selection}) arises from the fact that correlation functions do not have an analytic form, and hence a sufficiently long realization is required to numerically approximate them accurately. 
Nevertheless, since the implementation of Colored-nLIM allows for parallelization over $\gammanLIM$, a $\gamma$-selection remains feasible, provided sufficient computational resources are available.

The nLIM methodology can, in principle, be extended to incorporate polynomial dynamics of any order. 
However, this involves much higher-order correlation functions to estimate the dynamics. 
In practice, accurately obtaining these correlation functions may require substantially more observational data, which could be challenging. 
Thus, a trade-off must be carefully considered between the polynomial order (i.e., the number of estimated parameters) and the available sampling.



\section*{Acknowledgments}
The authors would like to thank Yan-Ning Kuo for the valuable discussion that inspired this study.




\section*{Funding}
This study is supported by the Council for Science, Technology and Innovation (CSTI), Cross-ministerial Strategic Innovation Promotion Program (SIP), the 3rd period of SIP "Smart Energy Management System" (Grant Number JPJ012207, funding agency: JST).


\appendix

\section{ Proof of nLIMs } 

In this appendix, we give a proof of Colored-nLIM (i.e., Eqs. (\ref{Eq:Colored-nLIM-E}) to (\ref{Eq:Colored-nLIM-FDR})). 
The White-nLIM follows the same workflow and is easier, so we leave it to the readers.

The Fokker-Planck equation of Eq. (\ref{Eq:Colored-nLIM}) reads
\begin{align*}
    \frac{\partial}{\partial t}P 
    &= -\sum_i \frac{\partial}{\partial x_i} F_i(x) P - \sum_{i,\,j} (\sqrt{2\bQ})_{ij} \frac{\partial}{\partial x_i} \eta_j P + \frac{1}{\bmgam}\sum_i \frac{\partial}{\partial \eta_i}\eta_i P + \frac{1}{2\bmgam^2}\sum_i \frac{\partial^2}{ \partial {\eta_i}^2} P \\
    &= \LFP P,
\end{align*}
where $\LFP$ is the Fokker-Planck operator and 
\[ F_i(x) = \sum_{j,\,k}\bB_{ijk}x_j x_k+\sum_j\bA_{ij}x_j +\bC_i.\]
Let $\pst$ denote a stationary distribution (i.e., $\LFP\pst=0$) and $\LFP^*$ denote the adjoint Fokker-Planck operator given by
\[ \LFP^* = \sum_i F_i \frac{\partial}{\partial x_i} + \sum_{i,\,j}(\sqrt{2\bQ})_{ij}\eta_j \frac{\partial}{\partial x_i} - \frac{1}{\bmgam}\sum_i \eta_i \frac{\partial}{\partial \eta_i} + \frac{1}{2\bmgam^2}\sum_i \frac{\partial^2}{ \partial {\eta_i}^2} .\]
Then, in the steady state, following \cite{Jung1985,Risken1989}, we have 
\begin{align} \label{Eq:Appendix-1}
    \langle f\left(\bx(\cdot +\tau),\bmeta(\cdot +\tau)\right),g\left(\bx(\cdot),\bmeta(\cdot)\right) \rangle &= \int_{\R^{2n}} f(x,\eta) e^{\tau \LFP}\pst(x,\eta) g(x,\eta) \, dxd\eta \nonumber \\
    &= \int_{\R^{2n}} g(x,\eta)\pst(x,\eta) e^{\tau \LFP^*} f(x,\eta) \, dxd\eta \nonumber \\
    &= \int_{\R^{2n}} g(x,\eta)\pst(x,\eta) \left(\bI + \tau\LFP^* + \frac{\tau^2}{2}\LFP^*\LFP^* \right) f(x,\eta) \,dxd\eta + O(\tau^3).
\end{align}

A straightforward computation shows that 
\begin{align} \label{Eq:Appendix-2}
    \LFP^* \, x_p &= \sum_{j,\,k}\bB_{pjk}x_j x_k+\sum_j\bA_{pj}x_j +\bC_p +\sum_j \eta_j \left(\sqrt{2\bQ}\right)_{pj}.
\end{align}
Therefore, by Eqs. (\ref{Eq:Appendix-1}) and (\ref{Eq:Appendix-2}), we compute the first-order derivative of the ($p$, $q$)-entry of $\bK$ as follows:
\begin{align*}
    \bK_{pq}'(0) &= \int_{\R^{2n}} x_q \pst(x,\eta) \LFP^* \, x_p \, dxd\eta \\
     &= \int_{\R^{2n}} x_q \pst(x,\eta) \left( \sum_{j,\,k}\bB_{pjk}x_j x_k+\sum_j\bA_{pj}x_j +\bC_p +\sum_j \eta_j \left(\sqrt{2\bQ}\right)_{pj} \right)\, dxd\eta \\
     &= \sum_{j,\,k}\bB_{pjk} \bM_{jkq}(0) +\sum_j\bA_{pj}\bK_{jq}(0) + \bC_p\bE_q + \sum_j \left(\sqrt{2\bQ}\right)_{pj}\left(\bK_{\bmeta\bx}\right)_{jq}.
\end{align*}
This is equivalent to Eq. (\ref{Eq:Colored-nLIM-K'}) in tensor form. 

The computation of $\bK''(0)$, $\bM'(0)$, and $\bM''(0)$ follows a similar procedure. 
We leave these derivations as an exercise to the readers.

Finally, with the fact that $\langle \bmeta(\cdot) \rangle = 0$, Eqs. (\ref{Eq:Colored-nLIM-E}), (\ref{Eq:Colored-nLIM-Q}) and (\ref{Eq:Colored-nLIM-FDR}) can be obtained by considering 
\[ 0 = \int_{\R^{2n}} x_j \LFP\pst\,dxd\eta = \int_{\R^{2n}} \pst \LFP^* \, x_j\,dxd\eta, \]
\[ 0 = \int_{\R^{2n}} \eta_i x_j \LFP\pst\,dxd\eta = \int_{\R^{2n}} \pst \LFP^* \,\eta_i x_j\,dxd\eta, \]
and 
\[ 0 = \int_{\R^{2n}} x_i x_j \LFP\pst\,dxd\eta = \int_{\R^{2n}} \pst \LFP^* \,x_i x_j\,dxd\eta, \]
respectively. They are also left to the readers as an exercise.

\bibliographystyle{abbrv}
\bibliography{Reference}

\begin{thebibliography}{10}

\bibitem{Aiken2013}
C.~M. Aiken, A.~Santoso, S.~McGregor, and M.~H. England.
\newblock {The 1970's shift in ENSO dynamics: A linear inverse model perspective}.
\newblock {\em Geophysical Research Letters}, 40(8):1612--1617, 2013.

\bibitem{Alexander2008}
M.~A. Alexander, L.~Matrosova, C.~Penland, J.~D. Scott, and P.~Chang.
\newblock Forecasting pacific ssts: Linear inverse model predictions of the pdo.
\newblock {\em Journal of Climate}, 21(2):385 -- 402, 2008.

\bibitem{Burgers2005}
G.~Burgers, F.-F. Jin, and G.~J. van Oldenborgh.
\newblock {The simplest ENSO recharge oscillator}.
\newblock {\em Geophysical Research Letters}, 32(13), 2005.

\bibitem{DelSole1996}
T.~DelSole.
\newblock Can quasigeostrophic turbulence be modeled stochastically?
\newblock {\em Journal of Atmospheric Sciences}, 53(11):1617 -- 1633, 1996.

\bibitem{DelSole2000}
T.~DelSole.
\newblock A fundamental limitation of markov models.
\newblock {\em Journal of the Atmospheric Sciences}, 57(13):2158 -- 2168, 2000.

\bibitem{Lorenzo2015}
E.~Di~Lorenzo, G.~Liguori, N.~Schneider, J.~C. Furtado, B.~T. Anderson, and M.~A. Alexander.
\newblock {ENSO} and meridional modes: A null hypothesis for pacific climate variability.
\newblock {\em Geophysical Research Letters}, 42(21):9440--9448, 2015.

\bibitem{Ham2019}
Y.-G. Ham, J.-H. Kim, and J.-J. Luo.
\newblock {Deep learning for multi-year ENSO forecasts}.
\newblock {\em Nature}, 573(7775):568--572, 2019.

\bibitem{Jin2007}
F.-F. Jin, L.~Lin, A.~Timmermann, and J.~Zhao.
\newblock Ensemble-mean dynamics of the {ENSO} recharge oscillator under state-dependent stochastic forcing.
\newblock {\em Geophysical Research Letters}, 34(3), 2007.

\bibitem{Jung1985}
P.~Jung and H.~Risken.
\newblock Motion in a double-well potential with additive colored gaussian noise.
\newblock {\em Zeitschrift f{\"u}r Physik B Condensed Matter}, 61(3):367--379, Sep 1985.

\bibitem{Kido2023}
S.~Kido, I.~Richter, T.~Tozuka, and P.~Chang.
\newblock Understanding the interplay between {ENSO} and related tropical {SST} variability using linear inverse models.
\newblock {\em Climate Dynamics}, 61(3):1029--1048, Aug 2023.

\bibitem{LeVeque2007}
R.~J. LeVeque.
\newblock {\em {Finite Difference Methods for Ordinary and Partial Differential Equations}}.
\newblock Society for Industrial and Applied Mathematics, 2007.

\bibitem{Lien2024}
J.~Lien, Y.-N. Kuo, and H.~Ando.
\newblock A linear inverse model for colored-gaussian noise, 2024.

\bibitem{Lorenz63}
E.~N. Lorenz.
\newblock Deterministic nonperiodic flow.
\newblock {\em Journal of Atmospheric Sciences}, 20(2):130 -- 141, 1963.

\bibitem{Lou2020}
J.~Lou, T.~J. O’Kane, and N.~J. Holbrook.
\newblock {A Linear Inverse Model of Tropical and South Pacific Seasonal Predictability}.
\newblock {\em Journal of Climate}, 33(11):4537 -- 4554, 2020.

\bibitem{Martinez2017}
C.~Martinez‐Villalobos, D.~J. Vimont, C.~Penland, M.~Newman, and J.~D. Neelin.
\newblock Calculating state-dependent noise in a linear inverse model framework.
\newblock {\em Journal of the Atmospheric Sciences}, 75:479--496, 2017.

\bibitem{Newman2011}
M.~Newman, S.-I. Shin, and M.~A. Alexander.
\newblock Natural variation in enso flavors.
\newblock {\em Geophysical Research Letters}, 38(14), 2011.

\bibitem{Penland1989}
C.~Penland.
\newblock Random forcing and forecasting using principal oscillation pattern analysis.
\newblock {\em Monthly Weather Review}, 117(10):2165--2185, 1989.

\bibitem{Penland1996}
C.~Penland.
\newblock A stochastic model of {I}ndo{P}acific sea surface temperature anomalies.
\newblock {\em Physica D: Nonlinear Phenomena}, 98(2):534--558, 1996.
\newblock Nonlinear Phenomena in Ocean Dynamics.

\bibitem{Penland1993}
C.~Penland and T.~Magorian.
\newblock { Prediction of Niño 3 Sea Surface Temperatures Using Linear Inverse Modeling }.
\newblock {\em Journal of Climate}, 6(6):1067--1076, 1993.

\bibitem{Penland1994}
C.~Penland and L.~Matrosova.
\newblock A balance condition for stochastic numerical models with application to the el niño-southern oscillation.
\newblock {\em Journal of Climate}, 7(9):1352--1372, 1994.

\bibitem{Penland1995}
C.~Penland and P.~D. Sardeshmukh.
\newblock The optimal growth of tropical sea surface temperature anomalies.
\newblock {\em Journal of Climate}, 8(8):1999--2024, 1995.

\bibitem{Risken1989}
H.~Risken.
\newblock {\em The {F}okker-{P}lanck equation}, volume~18 of {\em Springer Series in Synergetics}.
\newblock Springer-Verlag, Berlin, second edition, 1989.
\newblock Methods of solution and applications.

\bibitem{Sandri1996}
M.~Sandri.
\newblock {Numerical calculation of Lyapunov exponents}.
\newblock {\em Mathematica Journal}, 6(3):78--84, 1996.

\bibitem{Shin2021}
S.-I. Shin, P.~D. Sardeshmukh, M.~Newman, C.~Penland, and M.~A. Alexander.
\newblock Impact of annual cycle on {ENSO} variability and predictability.
\newblock {\em Journal of Climate}, 34(1):171--193, 2021.

\bibitem{sarkka_solin_2019}
S.~Särkkä and A.~Solin.
\newblock {\em Applied Stochastic Differential Equations}.
\newblock Institute of Mathematical Statistics Textbooks. Cambridge University Press, 2019.

\bibitem{Wang2017}
C.~Wang, C.~Deser, J.-Y. Yu, P.~DiNezio, and A.~Clement.
\newblock {El Ni{\~n}o and southern oscillation (ENSO): a review}.
\newblock {\em Coral reefs of the eastern tropical Pacific: Persistence and loss in a dynamic environment}, pages 85--106, 2017.

\bibitem{Wen1999}
C.~Wen, A.~Kumar, Y.~Xue, and M.~J. McPhaden.
\newblock {Changes in Tropical Pacific Thermocline Depth and Their Relationship to ENSO after 1999}.
\newblock {\em Journal of Climate}, 27(19):7230 -- 7249, 2014.

\bibitem{Zhao2024}
S.~Zhao, F.-F. Jin, M.~F. Stuecker, P.~R. Thompson, J.-S. Kug, M.~J. McPhaden, M.~A. Cane, A.~T. Wittenberg, and W.~Cai.
\newblock {Explainable El Ni{\~n}o predictability from climate mode interactions}.
\newblock {\em Nature}, 630(8018):891--898, 2024.

\end{thebibliography}

\end{document}